\documentclass[11pt]{amsart}
\usepackage[T1]{fontenc}
\usepackage[margin=1.15in]{geometry}
\usepackage{amssymb,mathrsfs,aliascnt}
\usepackage[colorlinks=true,linkcolor=blue,citecolor=blue,urlcolor=blue]{hyperref}
\usepackage[capitalise,nameinlink]{cleveref}

\newcommand{\sharedtheorem}[2]{%
  \newaliascnt{#1}{theorem}%
  \newtheorem{#1}[#1]{#2}%
  \aliascntresetthe{#1}%
  \crefname{#1}{#2}{#2s}}
\theoremstyle{plain}
\newtheorem{theorem}{Theorem}[section]
\crefname{theorem}{Theorem}{Theorems}
\sharedtheorem{proposition}{Proposition}
\sharedtheorem{lemma}{Lemma}
\sharedtheorem{corollary}{Corollary}
\theoremstyle{definition}
\sharedtheorem{definition}{Definition}
\sharedtheorem{remark}{Remark}

\DeclareMathOperator{\tr}{tr}
\DeclareMathOperator{\adj}{adj}
\DeclareMathOperator{\diag}{diag}
\DeclareMathOperator{\Hom}{Hom}
\DeclareMathOperator{\End}{End}
\DeclareMathOperator{\SL}{SL}
\DeclareMathOperator{\PSL}{PSL}
\DeclareMathOperator{\PGL}{PGL}
\DeclareMathOperator{\Lie}{Lie}
\newcommand{\Qbar}{\overline{\mathbb Q}}
\newcommand{\Fr}{\mathrm{Fr}}
\newcommand{\haff}{h_{\mathrm{aff}}}

\title{The Zilber--Pink conjecture for curves in \(Y(1)^3\)}
\author{Georgios Papas}
\address{Department of Mathematics, Weizmann Institute of Science, Israel}
\email{georgios.papas@weizmann.ac.il}
\date{18 November 2026}

\begin{document}

\begin{abstract}
We prove the Zilber--Pink conjecture for curves in \(Y(1)^3\). By work of
Habegger--Pila and Pila, it suffices to show that a curve defined over
\(\Qbar\), not contained in a proper special subvariety, contains only
finitely many points at which two isogeny relations with a common
coordinate hold. We prove a polynomial height bound for these points.
The proof uses \(G\)-functions centered at a smooth fiber and the
integral traces of the isogenies at places of proximity. In place of
the Andr\'e--Bombieri Hasse principle we use an interpolated version, in
which each class of places of proximity carries its own system of
equations.
\end{abstract}

\maketitle

\section{Introduction}\label{sec:height}

Let \(C\subset Y(1)^3\) be an irreducible curve that is Hodge generic,
that is, not contained in any proper special subvariety. Special
subvarieties of \(Y(1)^3\) are defined by conditions of the form
\(x_i=x\), with \(x\) a singular modulus, and \(\Phi_N(x_i,x_j)=0\), where
\(\Phi_N\) is the \(N\)-th modular polynomial; see \cite[\S2.1]{HP12}.
The Zilber--Pink conjecture \cite{Zil02,Pin05} predicts that only
finitely many points of \(C\) lie on special subvarieties of codimension
at least two; see \cite{Pil22} for background. Our main
result confirms this.

\begin{theorem}\label{thm:main}
Let \(C\subset Y(1)^3\) be an irreducible curve over \(\mathbb C\) that
is not contained in any proper special subvariety. Then only finitely
many points of \(C\) lie on special subvarieties of codimension at least
two. In other words, the Zilber--Pink conjecture holds for curves in
\(Y(1)^3\).
\end{theorem}

For curves not defined over \(\Qbar\), \Cref{thm:main} is due to Pila
\cite[Theorem~1.4]{Pil17}. For curves defined over \(\Qbar\), Habegger
and Pila \cite[Theorem~2]{HP12} proved finiteness for the special
subvarieties that are not strongly special, that is, those defined using
a singular modulus. As observed by Daw and Orr \cite[\S1.A]{DO25}, it
therefore suffices to consider the strongly special subvarieties. For
curves in \(Y(1)^3\), the points that remain are those at which two
isogeny relations with a common coordinate hold. We prove their
finiteness in \Cref{thm:two-relations}.

Habegger and Pila proved this finiteness for asymmetric curves
\cite[Theorem~1]{HP12}. The other known cases rest on the
\(G\)-functions method of Andr\'e \cite{Andre} and Bombieri
\cite{Bom81}, which was first used in unlikely intersections by Daw and
Orr \cite{DO21}; see \cite{DGS94} for background on \(G\)-functions and
\cite{And25} for a survey. Daw and Orr \cite{DO25} proved the
Zilber--Pink conjecture for curves whose closure in \((\mathbb P^1)^3\)
contains \((\infty,\infty,\infty)\). In earlier work \cite{Pap24}, the
author extended this to curves whose closure meets the boundary in a
point all of whose coordinates are \(\infty\) or singular moduli. In
joint work with Daw and Orr \cite{DOP25}, the author treated curves
whose closure meets a strongly special curve in the boundary.

In \cite{DO21,DO25,Pap24}, as in Andr\'e's original work
\cite[Chapter~X]{Andre}, the \(G\)-functions are centered at a point
where the families of elliptic curves or abelian varieties degenerate,
that is, at a point where the curve meets the boundary of a
compactification of the ambient Shimura variety. Such a center imposes
conditions on the points of interest at every finite place. These
conditions are used either to rule out \(v\)-adic proximity to the
center, as in \cite[Chapter~X]{Andre}, or to construct relations that
hold at all finite places simultaneously, as Daw and Orr do with the
Tate uniformization \cite[\S1.B]{DO25}. The idea to center the
\(G\)-functions instead at a point in the interior of the moduli space
goes back to Andr\'e \cite{And95} and Beukers \cite{Beu93}, who studied
the values of \(G\)-functions at CM points of \(Y(1)\). In particular,
relations among non-archimedean values of \(G\)-functions were used in
\cite[Th\'eor\`eme~1]{And95}.

In joint work with Daw and Orr \cite{DOP25}, the author took up this
point of view for curves in \(Y(1)^3\). There the \(G\)-functions are
centered at a point \(s_0\) of the curve whose coordinates are
isogenous, that is, at a smooth fiber, and relations among their values
at a point \(s\) are constructed from the reductions of the isogenies at
\(s\) and \(s_0\). At places of ordinary reduction, the two isogenies
give a single relation valid at all such places
\cite[Proposition~3.6]{DOP25}. At places of supersingular reduction, the
relations depend on the place. This gives a height bound for the points
with few places of supersingular proximity to \(s_0\)
\cite[Theorem~1.8]{DOP25}. The same point of view is developed for
families of abelian varieties in \cite{Pap25}.

In this paper we overcome the obstruction coming from supersingular
proximity by replacing the Andr\'e--Bombieri Hasse principle by its
interpolated version \cite[Theorem~5.17]{QM}, which corresponds to
\cite[Theorem~3.1]{OAI}. The Hasse principle bounds the height of points
at which a non-trivial global relation holds \cite[Chapter~VII,
\S5.2, Theorem]{Andre}. Its interpolated version descends from the
method of Bombieri \cite[\S10, Main Theorem]{Bom81} and Andr\'e
\cite[Chapter~VII, \S3.5, Proposition]{Andre}. In place of a global
relation, \cite[Theorem~5.17]{QM} allows the places of proximity to be
divided into classes, each with its own system of equations. In
return, each system must be shown to be non-trivial in the sense of
\cite[Definition~5.14]{QM}. Since the two isogenies have a common
source, their reductions give endomorphisms of the same elliptic curve,
and their commutator supplies one of these systems.

\subsection{Two isogeny relations}\label{subsec:setup}
From now on let \(C\) be defined over a number field \(K\), and write
\(j_1,j_2,j_3\) for its coordinates. After permuting the coordinates, two
isogeny relations with a common coordinate are \(\Phi_M(j_1,j_2)=0\) and
\(\Phi_N(j_1,j_3)=0\). We denote the first pair of indices by
\(e:=(1,2)\) and the second by \(f:=(1,3)\). A statement made for
\(k\in\{e,f\}\) is made for both pairs, and we then write \(k=(1,j)\),
with \(j\in\{2,3\}\).

Fix a smooth finite cover \(S\) of \(C\) over \(K\) carrying elliptic
schemes \(E_i\to S\) with \(j(E_i)=j_i\) \((i=1,2,3)\), and write
\(E_{i,s}\) for the fiber at \(s\in S(\Qbar)\). Let \(h_S\) be a Weil
height on \(S\) associated with an ample divisor \(H_S\) on the smooth
projective completion \(\overline S\). At a point where both
relations hold, the three fibers are isogenous, so they are either all
CM or all non-CM. In the first case the point is special. Such points are finite in
number by the Andr\'e--Oort conjecture for \(Y(1)^3\), proved by Pila
\cite{Pil11} after the case of curves in \(Y(1)^2\) due to Andr\'e
\cite{And98}; see also \cite[Theorem~2]{HP12}. We therefore consider
\[
 \mathcal S:=\left\{s\in S(\Qbar):
 \begin{array}{l}
 E_{1,s}\text{ is isogenous to }E_{2,s}\text{ and to }E_{3,s},\\
 E_{1,s}\text{ is not CM}
 \end{array}\right\}.
\]
To study this set we may assume that it is nonempty, since otherwise
there is nothing to prove. We may thus fix \(s_0\in\mathcal S\), which
we call the \emph{center}. Our main technical result is the following
height bound.

\begin{theorem}\label{thm:height}
Fix \(s_0\in\mathcal S\). There exist constants \(c_1,c_2>0\), depending
only on \(C\), \(S\) and \(s_0\), such that
\[
 h_S(s)\le c_1\,[K(s):K]^{c_2}
\]
for every \(s\in\mathcal S\).
\end{theorem}

Combined with the isogeny estimates of Masser and W\"ustholz
\cite{MW90}, in the effective form \cite[Theorem~1.4]{GR14}, as in the
proof of \cite[Proposition~5.15]{DO25},
\Cref{thm:height} gives a lower bound for Galois orbits. This bound
replaces \cite[Lemma~4.2]{HP12} in the counting argument of
\cite[Proposition~5.1]{HP12}, as in \cite[\S5]{DOP25}. That argument follows the strategy of Pila and Zannier \cite{PZ08} and
uses the counting theorem of Pila and Wilkie \cite{PW06}.
It gives the following; see \Cref{subsec:proofs}.

\begin{theorem}\label{thm:two-relations}
Let \(C\subset Y(1)^3\) be an irreducible curve defined over \(\Qbar\)
that is not contained in any proper special subvariety. Then \(C\)
contains only finitely many points \((j_1,j_2,j_3)\) for which there
exist \(M,N\ge1\) with
\[
 \Phi_M(j_1,j_2)=\Phi_N(j_1,j_3)=0.
\]
\end{theorem}

After permuting the coordinates, \Cref{thm:two-relations} is
\cite[Conjecture~1.2]{DOP25}. Together with the reduction above, it
proves \Cref{thm:main}; see \Cref{subsec:proofs}. The same
finiteness holds for points of curves in \(Y(1)^n\) satisfying two
isogeny relations with a common coordinate (\Cref{rem:three-coordinates}).

\subsection{Outline of the proof}\label{subsec:outline}
Let \(s\in\mathcal S\). We consider the embeddings of a number field
containing \(K(s)\) at which \(s\) is close to \(s_0\), in the sense of
\Cref{def:close}.

\subsubsection{Integral traces}
At each such embedding, parallel transport to the center turns the
moving isogeny of a pair \(k\in\{e,f\}\), composed with the dual of the
central isogeny, into an endomorphism of the homology of a central
fiber. Its trace is the value at \(s\) of a power series \(F_k\) built
from the \(G\)-functions; see
\eqref{eq:trace-series}. At finite places of good reduction it is the
trace of an endomorphism of the common reduction. At archimedean places
it is the trace of an endomorphism of the integral Betti lattice, and at
places of multiplicative reduction it is computed from Tate
uniformizations (\Cref{lem:toric-trace}). In every case it is an integer
of controlled size.

\subsubsection{Selection of places}
The interpolated Hasse principle requires a system of equations common to a
set of places carrying a proportion \(\gg1/[K(s):K]\) of the proximity
of \(s\) to \(s_0\). At archimedean places and above a fixed finite set
of primes, we group the places by the pair of integral traces. At the
remaining places we use the commutator of the two reduced endomorphisms
when it vanishes, as in \cite[Proposition~3.6]{DOP25}, and a vanishing
trace when one occurs. The places that are left are grouped by the
restriction of a local Frobenius lift, and paired with their conjugates
as in \cite[Lemma~5.1]{QM}. There are \(O([K(s):K])\) classes, so one of
them carries enough proximity (\Cref{prop:proximity}).

\subsubsection{Non-triviality}
It remains to show that the selected system is non-trivial. For the
systems on a single chart, this follows from the monodromy of the
elliptic families (\Cref{sec:single}). For a paired system, failure of
non-triviality forces the pair \((s,s^\tau)\) onto a coordinatewise
Hecke correspondence of bounded degree (\Cref{sec:hecke}). Frobenius then conjugates each of the two reduced endomorphisms to plus
or minus itself. Consequently they either commute, or one of them has
trace zero (\Cref{sec:frobenius}). Every place of the class would then
have been assigned to an earlier class, which is a contradiction.

\subsection{Notation}\label{subsec:conventions}
For \(s\in S(\Qbar)\) we set
\[
 d(s):=\max\bigl(2,[K(s):K]\bigr),\qquad
 h(s):=\max\bigl(2,h_S(s)\bigr),
\]
and write \(d=d(s)\) and \(h=h(s)\) when \(s\) is fixed. The symbols
\(\ll\) and \(O(\cdot)\), and the constants \(c,C\), depend only on
\(C\), \(S\), the center \(s_0\), and the fixed choices of
models, bases, discs and of the set \(\Sigma\) of \Cref{subsec:sigma};
they are independent of \(s\). We write \(h_{\mathbb P^1}\) for the
absolute logarithmic Weil height on \(\mathbb P^1(\Qbar)\), and
\(\haff\) for the affine height of \cite[Definition~5.10]{QM}.

For a place \(w\) of a number field \(L\), we write
\(\iota_w:L\hookrightarrow\mathbb C_w\) for the associated embedding,
where \(\mathbb C_w\) is \(\mathbb C\) or \(\mathbb C_p\) according as \(w\)
is archimedean or lies above \(p\). For \(y\in L[[x]]\), we write
\(\iota_w(y)\in\mathbb C_w[[x]]\) for the series obtained by applying
\(\iota_w\) to its coefficients. Let \(\mathcal Y\) be a finite family of
such power series. We write \(R_w(y)\) for the radius of convergence of
\(\iota_w(y)\), for \(y\in\mathcal Y\), and set
\[
 R_w(\mathcal Y):=\min_{y\in\mathcal Y}R_w(y).
\]
For a family of matrices of power series, \(\mathcal Y\) is the family
of their entries.

For a number field \(L\), we consider pairs \((q,\rho)\), where \(q\) is
a place of \(\mathbb Q\) and \(\rho:L\hookrightarrow\Omega_q\), with the
weights \(\omega_L(q,\rho)\) of \cite[\S5]{QM}. For a matrix \(A\) over
\(L\) we write \(\|\rho(A)\|_q\) for the maximum of \(|\rho(a)|_q\) over
the entries \(a\) of \(A\). For an elliptic
curve \(E\) over a number field and a finite place of good reduction,
\(\widetilde E\) denotes the geometric reduction, and
\(\Hom^0:=\Hom\otimes\mathbb Q\).

\subsection{Constant coordinates}\label{subsec:constant}

\begin{remark}\label{rem:constant-coordinates}
Suppose that a coordinate \(j_i\) is constant on \(C\). Then at most two
of the degrees of \(j_1,j_2,j_3\) on \(C\) are nonzero, so \(C\) is
asymmetric in the sense of \cite[\S1]{HP12}, and \Cref{thm:two-relations}
for \(C\) is a case of \cite[Theorem~1]{HP12}. \Cref{thm:height} also
holds in this case, with \(h_S(s)\ll1+\log d(s)\). Indeed, let
\(k\in\{e,f\}\) be a pair containing \(i\), with other index \(j\). Then
\(j_j\) is nonconstant, since otherwise \(C\) would lie in a modular
hypersurface \(\Phi_M(j_i,j_j)=0\). Fix an elliptic curve \(E_*\) over
\(K\) with \(j(E_*)=j_i\). After a fixed extension of \(K\), and an
extension of \(K(s)\) of bounded degree identifying the twist, the
isogeny estimate \cite{MW90}, \cite[Theorem~1.4]{GR14} with source
\(E_*\), as in the proof of \cite[Proposition~5.15]{DO25}, bounds the
minimal degree \(N\) of an isogeny \(E_*\to E_{j,s}\) by \(d(s)^c\). The
isogeny inequality \(h_F(E_{j,s})\le h_F(E_*)+\tfrac12\log N\) for the
stable Faltings height \cite[Corollary~2.1.4]{Ray85} and Silverman's
comparison of the Faltings and modular heights
\cite[Proposition~2.1]{Sil86} give
\(h_{\mathbb P^1}(j_j(s))\ll1+\log d(s)\). Since the pullback of
\(\mathcal O(1)\) under \(j_j:\overline S\to\mathbb P^1\) has positive
degree, the height machine \cite[Theorem~B.3.2]{HS00} gives
\(h_S(s)\ll1+h_{\mathbb P^1}(j_j(s))\).
\end{remark}

Henceforth every coordinate varies. Since \(C\) is Hodge generic, the
families \(E_1,E_2,E_3\) are then nonisotrivial and pairwise not
isogenous at the geometric generic point.

\section{Preliminaries}\label{sec:gfunctions}

We recall the \(G\)-functions centered at \(s_0\) and the interpretation
of their values in terms of periods, following \cite[\S\S2--3]{DOP25} and
\cite[\S3]{Pap25}.

\subsection{\texorpdfstring{\(G\)}{G}-functions centered at \texorpdfstring{\(s_0\)}{s0}}\label{subsec:bases}
We recall the construction of the \(G\)-functions associated with the
families \(E_i\) and the center \(s_0\); see \cite[\S2.6]{DOP25} and
\cite[\S3.3]{Pap25}. To prove \Cref{thm:height}, we may replace \(K\)
by a finite extension and \(S\) by a finite cover; see
\cite[\S2.2]{Pap24}. We may therefore assume that \(S\) carries the
covering datum of \cite[\S4.2]{DOP25}. It consists of a regular
projective model over the ring of integers of \(K\), semiabelian schemes
over it extending the \(E_i\), and a rational function \(x\) whose zeros
are simple and lie above \(s_0\). We fix one of these zeros as the lift
of \(s_0\), and denote it again by \(s_0\).

For each \(i\), we fix a Hodge basis \((\omega_i,\eta_i)\) of
\(H^1_{\mathrm{dR}}(E_i/S)\) on a fixed affine open neighborhood of
\(s_0\). That is, \(\omega_i\) is a section of the Hodge line
\(F^1H^1_{\mathrm{dR}}(E_i/S)\), and \(\langle\omega_i,\eta_i\rangle=1\)
for the cup product pairing. Such a basis exists, since \(s_0\) is an
interior point in characteristic zero, where the pairing is perfect.
With respect to this basis, the Gauss--Manin connection gives a
differential system
\begin{equation}\label{eq:gauss-manin}
 \frac{d}{dx}\,Y=A_i\,Y,\qquad A_i\in M_2(K(S)).
\end{equation}
Let \(Y_i\in M_2(K[[x]])\) be its normalized uniform solution, that is,
the solution with \(Y_i(0)=I_2\); see \cite[proof of
Theorem~3.4]{Pap25}. Its entries are \(G\)-functions
\cite[Lemma~2.3]{DOP25}. We call \(Y_i\) the matrix of \(G\)-functions
associated with \(E_i\), \(s_0\) and \((\omega_i,\eta_i)\). Since the cup
product is horizontal and \(\langle\omega_i,\eta_i\rangle=1\), we have
\(\det Y_i=1\).

\subsection{Places of proximity}\label{subsec:radii}
To ensure that analytic proximity to \(s_0\) implies equality of
reductions, we modify the bases as in \cite[Lemma~3.9]{Pap25},
simultaneously for the three families and with a single
\(\vartheta\in K^\times\). Set
\[
 H(x):=\frac{\vartheta}{\vartheta-x},
\]
and replace \((\omega_i,\eta_i)\) by \((H\omega_i,H^{-1}\eta_i)\). The
matrix \(\diag(H,H^{-1})\) preserves the symplectic form and equals the
identity at \(x=0\), and \(Y_i\) is replaced by \(\diag(H,H^{-1})\,Y_i\)
\cite[(11)]{Pap25}. This change does not affect the traces computed in
\Cref{sec:traces}, since the matrices of an isogeny at both of its
endpoints change accordingly. We keep the notation
\((\omega_i,\eta_i)\) and \(Y_i\) for the new data.

\begin{definition}[{\cite[Definition~3.12]{Pap25}}]\label{def:close}
Let \(s\in S(\Qbar)\) and let \(w\) be a place of a number field
containing \(K(s)\). Set
\[
 r_w:=\min\bigl\{1,R_w(Y_1,Y_2,Y_3)\bigr\}.
\]
We say that \(s\) is \emph{\(w\)-adically close to \(s_0\)} if
\(|x(s)|_w<r_w\) and \(s\) lies in the connected component through
\(s_0\) of the preimage under \(x\) of the open disc of radius \(r_w\)
centered at \(0\).
\end{definition}

By \cite[Lemma~3.11]{Pap25}, if \(s\) is \(w\)-adically close to \(s_0\)
at a finite place \(w\), then \(s\) and \(s_0\) have the same reduction
at \(w\). At archimedean places and at the finitely many finite places
above the set \(\Sigma\) of exceptional primes fixed in
\Cref{subsec:sigma}, we work instead on fixed closed inner discs around
\(s_0\). On these discs the local inverse of \(x\) through \(s_0\), the
comparison isomorphisms of \Cref{subsec:comparison} and the Tate
uniformizations are defined, and \(H^{\pm1}\) and the parallel
transports are bounded.

\subsection{\texorpdfstring{\(G\)}{G}-functions and periods}\label{subsec:comparison}
We recall the interpretation of the values of the \(Y_i\) in terms of
periods; see \cite[\S3.2]{DOP25} and \cite[\S3]{Pap25}. Let \(v\) be an
archimedean place of \(K\), or a finite place at which the central fibers
\(E_{i,s_0}\) have good reduction. Write \(H^1_v\) for Betti cohomology
with rational coefficients if \(v\) is archimedean, and for crystalline
cohomology of the reduction, tensored with the fraction field of the
Witt vectors of the residue field, if \(v\) is finite. For a fiber
\(E_{i,s}\) with good reduction at \(v\), the comparison isomorphism
between \(H^1_{\mathrm{dR}}(E_{i,s})\) and \(H^1_v(E_{i,s})\), after
extension of scalars, is due to Grothendieck at archimedean places and
to Berthelot and Ogus \cite{BO83} at finite places. In the basis
\((\omega_i,\eta_i)\) at \(s\) and a basis of \(H^1_v(E_{i,s})\), it is
recorded by the period matrix \(\Pi_v(E_{i,s})\) of
\cite[Definition~2.3]{Pap25}.

Let \(\Delta_v\) be the disc around \(s_0\) of \Cref{def:close}. At a
finite place, the fibers over \(\Delta_v\) have the same reduction as the
central fibers. By \cite[Lemmas~3.1 and~3.3]{DOP25}, the comparison
isomorphisms over \(\Delta_v\) are compatible with parallel transport
for the Gauss--Manin connection. We may thus identify
\(H^1_v(E_{i,s_0})\) with the horizontal sections of
\(H^1_{\mathrm{dR}}(E_i/S)\) over \(\Delta_v\), and use the resulting basis
of \(H^1_v(E_{i,s})\) for every \(s\in\Delta_v\). With these bases, the
relative period identity
\begin{equation}\label{eq:period-identity}
 \Pi_v(E_{i,s})=\iota_v\bigl(Y_i(x(s))\bigr)\,\Pi_v(E_{i,s_0})
 \qquad(s\in\Delta_v)
\end{equation}
holds \cite[Theorem~3.4(2)]{Pap25}. At archimedean places it follows
from the relative de Rham--Betti comparison, as first observed by
Andr\'e \cite{Andre}. At finite places it follows from the
de Rham--crystalline comparison of Berthelot and Ogus \cite{BO83} and
its version for families in Ogus's theory of convergent
\(F\)-isocrystals \cite[Theorems~3.8 and~3.10]{Ogu84}. Over primes
unramified in \(K/\mathbb Q\), this is due to Andr\'e \cite{And95}.

The comparison isomorphisms are functorial for isogenies. At a finite
place, an isogeny \(E_{i,s}\to E_{j,s}\) between fibers over \(\Delta_v\)
therefore acts on \(H^1_v\) through its reduction, which is an isogeny
between the reductions \(\widetilde E_i\) and \(\widetilde E_j\) of the
central fibers, as in \cite[\S3.2.1, (6)]{DOP25}. At places where the
central fibers have multiplicative reduction, crystalline cohomology is
replaced by the Tate uniformization; see \cite[\S3.2.2]{DOP25}.

\section{The two integral traces}\label{sec:traces}

\subsection{Transported isogenies}\label{subsec:trace-series}
Let \(k=(1,j)\in\{e,f\}\), and write
\[
 b_k:E_{1,s_0}\longrightarrow E_{j,s_0},\qquad
 a_k:E_{1,s}\longrightarrow E_{j,s}
\]
for the central and moving isogenies, of degrees \(N_{0,k}\) and
\(N_k\), respectively. In \Cref{subsec:degrees-isogenies} we take
\(a_k\) to be a generator of \(\Hom(E_{1,s},E_{j,s})\). On homology, in
the bases dual to the Hodge bases of \Cref{subsec:bases}, we write
\(\mathcal T_i(t)\) for parallel transport from the fiber of \(E_i\) at
the point with \(x=t\) to the central fiber \(E_{i,s_0}\)
\((i=1,2,3)\), \(A_k\) for the matrix of \(a_k\), and \(B_k^\vee\) for the
matrix of the dual central isogeny \(b_k^\vee\). For an independent
\(2\times2\) matrix variable \(P\), we set
\begin{equation}\label{eq:trace-series}
 M_k(t,P):=B_k^\vee\,\mathcal T_j(t)\,P\,\mathcal T_1(t)^{-1},\qquad
 F_k(t,P):=\tr M_k(t,P).
\end{equation}
Matrices act on column vectors, so that the matrix of a composite is
the product of the matrices. Write \(\Gamma_i\) for the column with
entries \(\omega_i,\eta_i\). Since \(Y_i\) is the normalized uniform
solution of the system \eqref{eq:gauss-manin} in this basis
\cite[proof of Theorem~3.4]{Pap25},
the column \(\Gamma_i^{\mathrm{hor}}:=Y_i^{-1}\Gamma_i\) is horizontal
wherever \(Y_i\) converges, and equals \(\Gamma_i\) at \(s_0\). A cycle
and its parallel transport have the same pairing with a horizontal
section. Hence, in the transport of the \(a\)-th dual basis vector at the
point with \(x=t\), the coefficient of the \(c\)-th dual basis vector at
\(s_0\) is \((Y_i(t)^{-1})_{ca}\). Since \(\det Y_i=1\), this gives
\[
 \mathcal T_i(t)=Y_i(t)^{-1}=\adj Y_i(t),\qquad
 \mathcal T_i(t)^{-1}=Y_i(t).
\]
In particular, the inverse introduces no denominator depending on
\(t\), and the entries of \(M_k(t,P)\) are polynomials in the entries of
\(P\), \(Y_1(t)\) and \(Y_j(t)\). Set \(\xi:=x(s)\). Then
\(M_k(\xi,A_k)\) is the endomorphism
\[
 H_1(E_{1,s_0})
 \xrightarrow{\ \mathcal T_1^{-1}\ } H_1(E_{1,s})
 \xrightarrow{\ a_k\ } H_1(E_{j,s})
 \xrightarrow{\ \mathcal T_j\ } H_1(E_{j,s_0})
 \xrightarrow{\ b_k^\vee\ } H_1(E_{1,s_0}).
\]
The next two lemmas compute its trace at an embedding at which \(s\) is
close to \(s_0\).

\begin{lemma}\label{lem:good-trace}
Let \(v\) be a place at which \(s\) is \(v\)-adically close to \(s_0\).
\begin{enumerate}
\item[(1)] If \(v\) is finite and the central fibers have good reduction
 at \(v\), then \(F_k(\xi,A_k)\in\mathbb Z\) and
 \(|F_k(\xi,A_k)|\le2\sqrt{N_{0,k}N_k}\).
\item[(2)] If \(v\) is archimedean and \(s\) lies in the fixed inner
 disc of \Cref{subsec:radii}, then \(F_k(\xi,A_k)\in\mathbb Z\) and
 \(|F_k(\xi,A_k)|\ll\sqrt{N_k}\).
\end{enumerate}
\end{lemma}

\begin{proof}
At a finite place of good reduction, the composite acts on the common
special elliptic curve \(\widetilde E_1\) as the endomorphism
\(\bar b_k^\vee\bar a_k\), by \Cref{subsec:comparison}. Its trace is
therefore an integer, also at supersingular places
\cite[Lemma~3.5 and (6)]{DOP25}; compare \cite[Lemma~4.9]{QM}. Since
this endomorphism has degree \(N_{0,k}N_k\), its trace has absolute
value at most \(2\sqrt{N_{0,k}N_k}\). At an archimedean place, parallel
transport preserves the integral Betti lattices, so the trace is again an
integer. On fixed compact discs around the center the transports are
bounded, while the Hodge norm of \(a_k\) is \(O(\sqrt{N_k})\).
\end{proof}

\subsection{Multiplicative places}\label{subsec:toric-trace}

\begin{lemma}\label{lem:toric-trace}
Suppose that the central fibers have multiplicative reduction at the
finite place \(v\), and that \(s\) is \(v\)-adically close to \(s_0\)
and lies in the fixed inner disc of \Cref{subsec:radii}. In compatible
split Tate uniformizations, let \(m_0\) and \(m\) denote
the signed torus exponents of \(b_k\) and \(a_k\), respectively. Then
\(m_0\mid N_{0,k}\), \(m\mid N_k\), and at this place
\begin{equation}\label{eq:toric-trace}
 F_k(\xi,A_k)=m\,\frac{N_{0,k}}{m_0}+m_0\,\frac{N_k}{m}\in\mathbb Z .
\end{equation}
\end{lemma}

\begin{proof}
After a fixed splitting extension of the local field and restriction to
a disc around \(s_0\), let \(\phi_i\) be the relative Tate uniformization of \(E_i\)
\cite[Chapter~V]{Sil94}, with coordinate \(w\) on the torus \((i=1,j)\). By
\cite[\S3.2.2, (8)]{DOP25}, the relation
\[
 \phi_i^*\gamma=\ell_i(\gamma)\,\frac{dw}{w}
\]
defines a horizontal quotient \(\ell_i\) of the relative de Rham
cohomology. Since the section \(dw/w\) of the target is horizontal,
parallel transport takes the generator of the dual line in homology at \(s\) to
the corresponding generator at \(s_0\). On these lines \(a_k\) acts by
\(m\), and \(b_k^\vee\) acts by \(N_{0,k}/m_0\), since
\(b_k^\vee b_k=[N_{0,k}]\). Hence \(M_k(\xi,A_k)\) has the eigenvalue
\[
 \lambda:=m\,\frac{N_{0,k}}{m_0}
\]
on the line through the generator at \(s_0\). Since the transports
preserve the pairing, we have
\[
 \det M_k(\xi,A_k)=N_{0,k}N_k,\qquad
 \tr M_k(\xi,A_k)=\lambda+\frac{N_{0,k}N_k}{\lambda},
\]
which is \eqref{eq:toric-trace}. The torus
exponents \(N_{0,k}/m_0\) and \(N_k/m\) of the dual isogenies are
integers, so \(m_0\mid N_{0,k}\) and \(m\mid N_k\), and the trace is an
integer. It does not change under extension of the local field.
\end{proof}

In particular \(|F_k(\xi,A_k)|\le 2N_{0,k}N_k\) in the usual absolute
value, so its logarithmic height is \(O(1+\log N_k)\). We apply
\Cref{lem:toric-trace} to both pairs \(e\) and \(f\), with the same
uniformization of the common source \(E_1\) and the same parallel
transport \(\mathcal T_1\).

\section{Degrees and de Rham coordinates}\label{sec:coordinates}

The degrees of the moving isogenies control the trace values. The
interpolated Hasse principle also requires the de Rham coordinates of the
moving isogenies to be small, globally and at the selected places. We
establish these bounds in this section.

\subsection{Degrees of the moving isogenies}\label{subsec:degrees-isogenies}
Let \(s\in\mathcal S\). Since the fibers are not CM, the group
\(\Hom(E_{1,s},E_{j,s})\) is free of rank one for each
\(k=(1,j)\in\{e,f\}\). A generator \(a_k\) is a primitive isogeny of
minimal degree, and its Galois conjugates over \(K(s)\) are \(\pm a_k\).
Consequently, after fixed base changes, the generators \(a_e\) and
\(a_f\) are defined over an extension of \(K(s)\) of degree at most four.
The elliptic isogeny estimate \cite{MW90}, \cite[Theorem~1.4]{GR14}, the
comparison between the Faltings and modular heights
\cite[Proposition~2.1]{Sil86}, and the height machine
\cite[Theorem~B.3.2]{HS00} give
\begin{equation}\label{eq:degree-bound}
 N_e,N_f\ll (dh)^c,
\end{equation}
as in the proof of \cite[Proposition~5.15]{DO25}. We fix a number field
\(L_s\) containing the fixed data, the point \(s\) and both moving
isogenies, with \([L_s:\mathbb Q]\ll d\).

\subsection{Coordinates on a multiplicative disc}\label{subsec:toric-coordinates}

\begin{lemma}\label{lem:toric-coordinates}
In the situation of \Cref{lem:toric-trace}, there is a fixed smaller
closed disc around \(s_0\) on which the matrix \(A_k\), in the original
algebraic bases, satisfies
\begin{equation}\label{eq:toric-coordinates}
 \log\max\bigl(1,\|A_k\|_v\bigr)\le C_v ,
\end{equation}
where \(C_v\) is independent of \(s\), of its local field of definition
and of \(N_k\).
\end{lemma}

\begin{proof}
For \(i\in\{1,j\}\), set
\[
 c_i:=\ell_i(\omega_i),\qquad c'_i:=\ell_i(\eta_i).
\]
Since the uniformization takes an invariant differential to a nonzero
multiple of \(dw/w\), we have \(c_i(s_0)\ne0\). On a sufficiently small
disc we may therefore set
\[
 \alpha_i:=c_i^{-1}\omega_i,\qquad
 \beta_i:=c_i\eta_i-c'_i\omega_i ,
\]
so that \(\ell_i(\alpha_i)=1\), \(\ell_i(\beta_i)=0\) and
\(\langle\alpha_i,\beta_i\rangle=1\). Thus \(\alpha_i\) spans the Hodge
line and \(\beta_i\) spans the kernel of \(\ell_i\). At \(s\), the isogeny
\(a_k\) preserves these two lines in the respective fibers. Its torus
exponent gives \(a_k^*\alpha_j=m\,\alpha_1\), and since
\(\langle a_k^*\alpha_j,a_k^*\beta_j\rangle=N_k\), we obtain
\(a_k^*\beta_j=(N_k/m)\beta_1\). Hence
\begin{equation}\label{eq:tate-diagonal}
 [a_k^*]_{(\alpha_j,\beta_j)\to(\alpha_1,\beta_1)}
   =\diag(m,N_k/m).
\end{equation}
Both entries are \(v\)-integral by \Cref{lem:toric-trace}, including when
the residue characteristic divides \(N_k\). The two changes of basis and
their inverses have bounded analytic entries on a fixed smaller closed
disc. This gives \eqref{eq:toric-coordinates} for \(a_k^*\), and, by
transposition in the dual bases, for \(A_k\).
\end{proof}

The bases in this proof can be chosen for the three factors on the same
disc.

\subsection{Global height of the two matrices}\label{subsec:elliptic-heights}
Let \(L\) be a number field containing the data.

\begin{lemma}\label{lem:elliptic-global}
For \(k=(1,j)\in\{e,f\}\), the matrix \(A_k\) of the moving isogeny, in the
fixed algebraic de Rham bases, satisfies
\begin{equation}\label{eq:elliptic-global}
 \haff(A_k)\le C h+\tfrac12\log N_k .
\end{equation}
In particular \(\haff(A_e,A_f)\ll h+\log d\).
\end{lemma}

\begin{proof}
We first fix charts. For each \(E_i\), choose a short Weierstrass
equation on an open neighborhood of \(s_0\). Repeating this around each
of the finitely many non-CM points omitted by this open gives finitely
many charts, and it suffices to work on one of them. Write \(X_0,Y_0\)
for its coordinates and \(c_4,\Delta\) for its invariants, and set
\(\omega_0:=dX_0/(2Y_0)\). Since \(s_0\) is smooth and non-CM, we may
assume that \(c_4\Delta\ne0\) on the open. The changes of basis between
\((\omega_0,X_0\omega_0)\) and the chosen algebraic basis, in both
directions, then have entries in \(K(S)\) that are regular on the open.

\emph{Finite places.} Fix \(s\) and a pair \((q,\rho)\) with \(q\)
finite. We pass to a finite extension of the local field over which
both fibers have good reduction or both have split multiplicative
reduction; the two reduction types agree, since the fibers are
isogenous. Let \(R\) be its valuation ring. In the good case we take a
smooth integral Weierstrass equation, and in the multiplicative case the
standard Tate equation. Write \(X,Y\) for its coordinates,
\(a_1,\ldots,a_6\) for its coefficients, and
\(\Omega:=dX/(2Y+a_1X+a_3)\) for its invariant differential. We claim
that pullback under \(a_k\) is integral with respect to the bases
\((\Omega,X\Omega)\) at the two endpoints.

In the good case, \((\Omega,X\Omega)\) is an \(R\)-basis of the de Rham
cohomology of the smooth proper model. Indeed, let \(\pi:=-X/Y\) be the
uniformizer at the identity. Then \(\Omega=(1+O(\pi))\,d\pi\), and
\(X\Omega\) has principal part \(\pi^{-2}d\pi\). The coefficient of
\(\pi^{-1}d\pi\) vanishes by the residue theorem on the generic fiber,
and all Laurent coefficients lie in \(R\). Subtracting \(d(-1/\pi)\)
near the identity therefore gives an integral de Rham cocycle. Its image
in \(H^1(\mathcal O)\) is represented by the principal part \(\pm1/\pi\),
which pairs with \(\Omega\) to a unit by the residue description of
Serre duality. Since \(\Omega\) spans the Hodge line, the Hodge exact
sequence shows that \((\Omega,X\Omega)\) is a basis. The isogeny extends
to the abelian schemes over \(R\), so its pullback is integral with
respect to these bases.

For a Tate fiber of parameter \(\mathfrak q\), the uniformization takes \(\Omega\)
to \(dw/w\). On the annulus \(|\mathfrak q|_v<|w|_v<1\), the coefficient of \(w^0\)
in the Laurent expansion of \(X(w,\mathfrak q)\) is
\[
 \kappa(\mathfrak q):=-2\sum_{n\ge1}\frac{n\mathfrak q^n}{1-\mathfrak q^n}\in R
\]
\cite[(1)--(4)]{Tat95}, so the annular residue of \(X(w,\mathfrak q)\,dw/w\) is
\(\kappa(\mathfrak q)\). Restriction to the annulus takes the generator \(dw/w\) of
the de Rham cohomology of the torus to a differential of residue one.
Hence \((X-\kappa(\mathfrak q))\Omega\) spans the kernel of the uniformization
quotient in characteristic zero. By the principal-part calculation at
the identity, its cup product with \(\Omega\) is the same sign \(\pm1\)
for every such equation. The proof of \eqref{eq:tate-diagonal}
therefore applies to the basis \((\Omega,(X-\kappa(\mathfrak q))\Omega)\) at each
endpoint. Since \(\kappa(\mathfrak q)\in R\), the change of basis from
\((\Omega,X\Omega)\) is integral with integral inverse. This proves the
claim at every multiplicative fiber, whether or not \(s\) is close to
the center.

It remains to bound the change of basis back to the short equation.
Write \(X_0=u^2X+r\) for the change of variables. Then
\[
 \omega_0=u^{-1}\Omega,\qquad
 X_0\omega_0=uX\Omega+(r/u)\Omega,\qquad
 12r=u^2b'_2,
\]
where \(b'_2:=a_1^2+4a_2\in R\); the last identity holds because the
short equation has \(b_2=0\). Consequently this change of basis and its
inverse have maximum norm at most
\(|12|_q^{-1}\max(|u|_q,|u|_q^{-1})\). In the good case
\(\Delta(s)=u^{12}\Delta'\) with \(|\Delta'|_q=1\), and in the Tate case
\(c_4(s)=u^4c'_4\) with \(|c'_4|_q=1\), since \(c'_4\in1+\mathfrak q\mathbb Z[[\mathfrak q]]\).
In either case
\[
 \bigl|\log|u|_q\bigr|
 \le \tfrac1{12}\bigl|\log|\rho(\Delta(s))|_q\bigr|
       +\tfrac14\bigl|\log|\rho(c_4(s))|_q\bigr|,
\]
where the right-hand side involves only \(\Delta(s),c_4(s)\in K(s)\). For
\(g\in\{\Delta,c_4\}\), the product formula and the height machine give
\[
 \sum_{\substack{(q,\rho)\\q\ne\infty}}\omega_L(q,\rho)\,\bigl|\log|\rho(g(s))|_q\bigr|
 \le 2h_{\mathbb P^1}(g(s))\ll h .
\]
For each entry \(g\) of a change of basis to the chosen algebraic basis,
or of its inverse, the sum of
\(\omega_L(q,\rho)\log\max(1,|\rho(g(s))|_q)\) over the pairs \((q,\rho)\) with \(q\) finite is at
most \(h_{\mathbb P^1}(g(s))\ll h\), including when \(g(s)=0\). Moreover
\(\sum_{(q,\rho):\,q\ne\infty}\omega_L(q,\rho)\log|12|_q^{-1}=\log12\). Applying the
changes of basis at both endpoints to the integral pullback matrix, we
obtain
\begin{equation}\label{eq:elliptic-finite}
 \sum_{\substack{(q,\rho)\\q\ne\infty}}\omega_L(q,\rho)\log\max\bigl(1,\|\rho(A_k)\|_q\bigr)
 \ll h .
\end{equation}

\emph{Archimedean places.} With respect to the polarized Hodge norms,
the adjoint of \(a_k^*\) is \((a_k^\vee)^*\), so the operator norm of
\(a_k^*\) is \(\sqrt{N_k}\). Let \(K_\rho(s)\ge1\) bound the product of
the norms of the chosen basis at the target and of its dual at the
source. Then
\[
 \|\rho(A_k)\|_\infty\le K_\rho(s)\sqrt{N_k}.
\]
The archimedean calculation in the proof of \cite[Lemma~6.9]{QM} applies
to the two polarized variations of weight one. Near a boundary point or a
pole of the bases, with local parameter \(t\), it gives
\[
 K_\rho\le C|t|^{-M}\bigl(1+|\log|t||\bigr)^M ,
\]
by the fixed pole orders of the bases and their duals in the canonical
extensions and the Hodge norm estimate there. On the remaining compact
set \(K_\rho\) is bounded. Thus
\(\log K_\rho(s)\ll1+\lambda_{\mathfrak D,\rho}(s)\) for a fixed
effective divisor \(\mathfrak D\) supported at the omitted points and a
nonnegative local height \(\lambda_{\mathfrak D,\rho}\), and the weighted
archimedean sum is \(O(h)\) by the height bound for \(\mathfrak D\).

Since the archimedean weights sum to one, this and
\eqref{eq:elliptic-finite} prove \eqref{eq:elliptic-global}. The bound
for the two matrices follows by adding their heights and using
\eqref{eq:degree-bound}.
\end{proof}

\subsection{Selected places}\label{subsec:selected-places}
At almost every finite place, the center, the bases and their inverses
extend over a smooth integral Weierstrass chart around the reduction of
\(s_0\). Points in its residue disc have good models and integral
pullbacks there, so that \(\|\rho(A_k)\|_q\le1\). At each of the
finitely many remaining places at which the central fibers have good
reduction, we fix an inner disc on which the functions giving the two
changes of basis to short equations, and the inverses of their
discriminants, are bounded. The good-reduction calculation in the proof
of \Cref{lem:elliptic-global} then bounds \(\|\rho(A_k)\|_q\) by a fixed
constant. At places where the central fibers have multiplicative
reduction we use \Cref{lem:toric-coordinates}. Finally, on closed
archimedean discs around the center \(K_\rho(s)\) is bounded, so the
local logarithm is \(O(1+\log N_k)\).

Let \(\mathcal V\) be a set of pairs \((q,\rho)\) at which \(\rho(s)\)
lies in these discs. Since the weights of the pairs in \(\mathcal V\)
above each rational place sum to at most one, \eqref{eq:degree-bound}
gives
\begin{equation}\label{eq:selected-matrices}
 \sum_{(q,\rho)\in\mathcal V}\omega_{L_s}(q,\rho)
   \log\max\bigl(1,\|\rho(A_e,A_f)\|_q\bigr)
 \ll 1+\log h+\log d .
\end{equation}
If \(\sigma_\rho\) is a local Frobenius lift and
\(\tau=\rho^{-1}\sigma_\rho\rho\), as in \Cref{subsec:partition}, then
\[
 \|\rho(\tau(A_k))\|_p=\|\sigma_\rho(\rho(A_k))\|_p=\|\rho(A_k)\|_p .
\]
We use this on the doubled chart.

\section{Common equations and selection}\label{sec:selection}

In this section we divide the places of proximity into classes, attach a
system of equations to each class, and verify hypotheses (i)--(iii) of
\cite[Theorem~5.17]{QM} for a class carrying enough proximity.

\subsection{Variables and charts}\label{subsec:variables}
We use the fixed affine chart of \cite[\S5]{QM}, defined on an affine
open subset \(S_{\mathrm{aff}}\subseteq S\) containing \(s_0\), on which
the bases of \Cref{subsec:bases} are also defined. Set
\[
 \mathcal S^\circ:=\{s\in\mathcal S\cap S_{\mathrm{aff}}(\Qbar):x(s)\ne0\}.
\]
The complement \(\mathcal S\setminus\mathcal S^\circ\) is finite and
independent of the moving point; it contains \(s_0\), which remains the
center of the series. Its points have bounded height, so it suffices to
prove \Cref{thm:height} for \(s\in\mathcal S^\circ\), after enlarging
\(c_1\). Let \(s\in\mathcal S^\circ\), so that \(\xi=x(s)\ne0\). Write \(z\) for its
coordinate functions,
\(Z\) for the corresponding independent variables, and \(\zeta(t)\) for
the inverse series through \(s_0\). Let \(P_e,P_f\) be independent
\(2\times2\) matrix variables.

The \emph{fixed series} are the entries of \(Y_1,Y_2,Y_3\) and the components of
\(\zeta(t)\). By \Cref{subsec:trace-series}, the entries of
\(\mathcal T_i^{\pm1}\) are among them, up to sign. Their coefficients
lie in a fixed number field, so each has finitely many conjugates. As in
\cite[Definition~5.11]{QM}, let \(\mathbf B\) be a tuple of independent
variables, one for each fixed series and each of its conjugates, and let
\(\mathbf B'\) be a second such tuple. Each series in the systems below
is obtained from a polynomial in \(t\), \(U\) and \(\mathbf B\) (on the
doubled chart, also \(\mathbf B'\)) by replacing each variable of
\(\mathbf B\) by the corresponding series in \(t\), and each variable of
\(\mathbf B'\) by the same series in \(t\mathcal R\). Such a polynomial
is a \emph{polynomial representative} of the series.

A \emph{single-chart system} has the variables
\[
 \eta:=(t,U),\qquad U:=(Z,P_e,P_f),
\]
whose value at \(s\) is \(\eta(s):=(\xi,z(s),A_e,A_f)\). It contains all
components of \(Z-\zeta(t)\).

A system on the \emph{doubled chart} depends on an embedding \(\tau\)
and has the variables
\[
 \eta:=(t,U),\qquad U:=(Z,P_e,P_f,Z^\tau,P'_e,P'_f,\mathcal R),
\]
whose value at \(s\) is
\[
 \eta(s):=\bigl(\xi,z(s),A_e,A_f,z^\tau(s^\tau),\tau(A_e),\tau(A_f),
 R(\xi,\tau)\bigr),\qquad\text{where}\quad R(\xi,\tau):=\frac{\tau(\xi)}{\xi}.
\]
It contains all components of \(Z-\zeta(t)\) and of
\(Z^\tau-\zeta^\tau(t\mathcal R)\).

In both cases the algebraic chart equations are those of \cite[\S5]{QM},
with \(x=t\) and, on the doubled chart, \(x^\tau=t\mathcal R\); the
endomorphism entries there are replaced by the entries of \(P_e,P_f\)
and, on the doubled chart, of \(P'_e,P'_f\).

\subsection{The exceptional primes}\label{subsec:sigma}
Let \(\Sigma\) be a fixed finite set of rational primes containing the
primes of bad reduction for the central models and for the charts, and
the primes at which two distinct zeros of \(x\) have the same reduction.
Thus, for \(p\notin\Sigma\), the residue disc through \(s_0\) contains no
other zero of \(x\). \(\Sigma\) also contains the primes of \Cref{lem:lifting}. These lists depend only on the fixed curves, their
finitely many conjugates and the center, so \(\Sigma\) is chosen before
any embedding is selected.

With this choice, the construction of \cite[Definition~6.19]{QM} and the
proofs of \cite[Lemmas~6.22--6.23]{QM} apply verbatim to the elliptic
families. They use only the following facts:
\begin{itemize}
\item the fixed series and their conjugates converge on \(|t|_p<1\) for
 \(p\notin\Sigma\), by Andr\'e's theorem for geometric differential
 equations \cite[Chapter~V, Appendix, pp.~110--111]{Andre}, after
 enlarging \(\Sigma\);
\item the inverse chart is integral on the residue disc;
\item the factorial cutoff and the product formula;
\item local Weil functions for the divisor \([s_0]\).
\end{itemize}

\subsection{The partition}\label{subsec:partition}
Let \(M_s\) be the normal closure of \(L_s\). We start from the
embeddings on the discs of \cite[Lemma~6.23]{QM}, that is, the initial
selection of \cite[Definition~6.19]{QM}, and, as in \cite[\S6]{QM}, we
retain all their extensions \(\rho:M_s\hookrightarrow\Omega_q\). By
\cite[Lemma~5.4]{QM} this preserves the weighted sums. We partition this
extended set; each embedding is assigned to the first of the following
classes that applies.

\begin{enumerate}
\item[(1)] \emph{Archimedean embeddings and embeddings above
\(p\in\Sigma\).} We group them by the ordered pair \(\mathbf
T=(T_e,T_f)\) of trace values \(T_k=F_k(\xi,A_k)\), which are integers by
\Cref{lem:good-trace,lem:toric-trace}. At the exceptional primes, these
lemmas are applied on the fixed inner discs of \Cref{subsec:radii}. The
equations are
\begin{equation}\label{eq:fixed}
 F_e(t,P_e)-T_e=0,\qquad F_f(t,P_f)-T_f=0.
\end{equation}
For a given \(s\) there are \(O(d)\) labels \(\mathbf T\), since \(L_s\)
has \(O(d)\) embeddings above the fixed finite set
\(\Sigma\cup\{\infty\}\). This count includes the multiplicative places.

\item[(2)] \emph{Good embeddings above \(p\notin\Sigma\) with commuting
reductions.} The selected disc identifies each moving reduction with the
central reduction \(\widetilde E_i\). For \(k\in\{e,f\}\), we set, in the
geometric endomorphism algebra of \(\widetilde E_1\),
\[
 D_k:=\bar b_k^{-1}\bar a_k,\qquad\text{so that}\qquad
 N_{0,k}D_k=\bar b_k^\vee\bar a_k;
\]
here and below, the inverse of an isogeny is taken in \(\Hom^0\). If
\(D_eD_f=D_fD_e\), we assign the embedding to the
\emph{commutator system}
\begin{equation}\label{eq:comm}
 \bigl[M_e(t,P_e),M_f(t,P_f)\bigr]=0,
\end{equation}
that is, to the four equations given by the entries of this
commutator. These equations hold at the embedding. Indeed, the
comparison for the common source \(E_1\) identifies the two evaluated
matrices, in the same basis, with the matrices of \(N_{0,e}D_e\) and
\(N_{0,f}D_f\), which commute. This is the commutation argument of
\cite[Proposition~3.6, (12)]{DOP25}. It applies at every ordinary place,
where the geometric endomorphism algebra is an imaginary quadratic
field \cite[\S\S III.9 and~V.3]{Sil09}, and at every supersingular place at which the two reductions
commute.

\item[(3)] \emph{Good embeddings at which a trace vanishes.} If
\(F_k(\xi,A_k)=0\) at the embedding for some \(k\in\{e,f\}\), we assign
it to the system
\begin{equation}\label{eq:zero}
 F_k(t,P_k)=0,
\end{equation}
taking \(k=e\) if both traces vanish. This gives at most two further classes.

\item[(4)] \emph{The remaining good embeddings.} Let
\(\rho:M_s\hookrightarrow\Qbar_p\) be such an embedding. We choose \(\sigma_\rho\in\operatorname{Gal}(\Qbar_p/\mathbb Q_p)\)
lifting the absolute Frobenius \(u\mapsto u^p\) of the residue field, and
set
\[
 \tau=\tau_\rho:=\rho^{-1}\sigma_\rho\rho|_{L_s}:
 L_s\hookrightarrow M_s .
\]
This conjugation by a Frobenius lift follows the proof of
\cite[Proposition~4.6]{OAI}, as developed in \cite[\S5.1]{QM}. We group
these embeddings by \(\tau\). The equations of the \emph{paired system}
are
\begin{equation}\label{eq:paired}
 F_k(t,P_k)-F_k^\tau(t\mathcal R,P'_k)=0
 \qquad(k=e,f),
\end{equation}
where \(F_k^\tau\) is obtained by applying \(\tau\) to every coefficient
of \(F_k\), including the fixed central data. The evaluation identity in
the proof of \cite[Lemma~5.1]{QM} transports the center to \(s_0^\tau\)
and fixes each integral trace. It gives \eqref{eq:paired} at
\(\eta(s)\) and \(|\rho(R(\xi,\tau))|_p=1\).
\end{enumerate}

Each class is used with its own equations and chart. The commutator
system \eqref{eq:comm} is a single system on a single chart, independent
of the prime, so it adds at most one class, and the systems
\eqref{eq:fixed} and \eqref{eq:zero} also use a single chart. The number
of classes of the paired system \eqref{eq:paired} is at most the number
\([L_s:\mathbb Q]\) of embeddings \(L_s\hookrightarrow M_s\), over all
rational primes together. In total there are \(O(d)\) classes. We apply
\cite[Theorem~5.17]{QM} with coefficient field \(E=M_s\) for every
class; its constants do not depend on \([E:\mathbb Q]\).

\subsection{Proximity}\label{subsec:proximity}

\begin{proposition}\label{prop:proximity}
There is a fixed polynomial \(Q_{\mathrm{thr}}\) with the following
property. If \(h\ge Q_{\mathrm{thr}}(d)\), then one of the classes of
\Cref{subsec:partition}, with set of embeddings \(\mathcal V\),
satisfies
\begin{equation}\label{eq:proximity}
 \sum_{(q,\rho)\in\mathcal V}\omega_{M_s}(q,\rho)\log|\rho(\xi)|_q^{-1}
 \gg h/d .
\end{equation}
\end{proposition}

\begin{proof}
After the fixed field extension making \(s_0\) rational, set
\(D_0:=[s_0]\) on \(\overline S\), and fix \(n_0\ge1\) with
\(\deg(n_0D_0-H_S)>0\). Since \(n_0D_0-H_S\) is ample, we have
\(h_{D_0}(s)\ge ch-C\). At a prime outside \(\Sigma\), the local equation
\(x\) generates the ideal of \(D_0\) on the residue disc through the
center. Summing its local heights as in \cite[Lemma~6.23]{QM}, we
therefore obtain \(h_{D_0}(s)+O(1)\) for the total proximity on the
preliminary discs.

At an exceptional place, replacing the disc by a fixed inner disc of
radius \(r\le1\) discards at most \(\log(1/r)\) times the total weight of
the embeddings there, hence \(O(1)\) in all. At the other places,
\cite[Lemma~6.22]{QM} bounds the loss from discarding the embeddings
with \(0<\log|\rho(\xi)|_p^{-1}\le 2\log(p)/(p-1)\) by
\(O(1+\log^2(dh))\). Thus, for \(h\) above a fixed polynomial in \(d\),
the retained sum is \(\gg h\). Dividing it among the \(O(d)\) classes
gives \eqref{eq:proximity}, as in the proof of \cite[Theorem~5.6]{QM},
which follows \cite[Lemma~4.5 and Proposition~4.6]{OAI}.
\end{proof}

\subsection{Hypotheses (i)--(iii) of the interpolated Hasse principle}\label{subsec:hypotheses}
Let \(\mathcal V\) be the class of \Cref{prop:proximity}, with its
system of equations and the value \(\eta(s)\).

For (i), the chart equations have fixed degrees and coefficient heights,
since only finitely many conjugates of the fixed data occur. The fixed
functions \(z\) and \(x\) have height \(O(h)\) at \(s\). Since
\(z^\tau(s^\tau)=\tau(z(s))\), the conjugate point has the same height.
On the doubled chart, moreover,
\[
 h_{\mathbb P^1}\bigl(R(\xi,\tau)\bigr)\le
 h_{\mathbb P^1}(\tau(\xi))+h_{\mathbb P^1}(\xi)=2h_{\mathbb P^1}(\xi).
\]
Since \Cref{lem:elliptic-global} applies to both matrices and their
conjugates, it follows that \(\haff(\eta(s))\ll h+\log d\).

For (ii), \Cref{prop:proximity} bounds the selected proximity from
below. The selected sum for the matrices is \(O(1+\log h+\log d)\) by
\eqref{eq:selected-matrices}. The fixed chart coordinates are integral
on the good discs outside \(\Sigma\) and bounded on the finitely many
other inner discs. In the paired system their conjugates have the same
local norms under the selected Frobenius lifts, and
\(|\rho(R(\xi,\tau))|_p=1\). This gives the same upper bound for the
selected sum of the full tuple \(U(s)\).

For (iii), the elliptic differential systems and inverse charts are of
the form required in \cite[Definition~5.11]{QM}. Since
\(\mathcal T_i^{-1}=\adj\mathcal T_i\), the entries of \eqref{eq:comm}
have degree at most four in the transport entries and are bilinear in
the two moving matrices, and all polynomial representatives have
bounded degree. The only scalar coefficients of \eqref{eq:fixed} that
vary with \(s\) are the two trace values, which have logarithmic height
\(O(1+\log(dh))\) by \Cref{lem:good-trace,lem:toric-trace} and
\eqref{eq:degree-bound}. The integral connection at almost all primes,
the algebraic inverse charts and the \(G\)-function estimates give the
bounds of \cite[Lemma~6.16]{QM}. At a good prime, a product of
coefficients of total order \(l\) has the same factor
\(\exp(l\log p/(p-1))\). Therefore the selected coordinate bounds and
the estimates of \cite[Propositions~6.25 and~6.28]{QM} give the
coefficient and tail bounds of (iii) for the additional systems as well;
for the commutator system see \Cref{rem:commutator-tails}.

\begin{remark}\label{rem:commutator-tails}
The proof of \cite[Proposition~6.28]{QM} writes each coefficient of
\(t^n\) as a fixed linear combination of the trace values and of Taylor
coefficients of products of the fixed series, and uses that the series
are linear in the matrix variables. This holds for \eqref{eq:fixed},
\eqref{eq:zero} and \eqref{eq:paired}. The commutator system
\eqref{eq:comm} is bilinear in \((P_e,P_f)\) and of degree four in
the transport entries. The same proof applies with the following
changes.
\begin{itemize}
\item The coefficient of \(t^n\) in a product of four fixed series is a
 sum of at most \((n+1)^3\) products of Taylor coefficients.
\item At \(p\notin\Sigma\), the ultrametric inequality and
 \cite[Lemma~6.16]{QM} give the same factor \(\exp(n\log p/(p-1))\) as
 for a single series.
\item At the finitely many other places, the factor
 \((n+1)^3\le e^{3n}\) is absorbed into the numbers \(b_v\) of
 \cite[Theorem~5.17(iii)]{QM}.
\item Set \(M_v(s):=\max(1,\|\rho(U(s))\|_q)\). The dependence on \(U\)
 is through \(M_v(s)^2\) instead of \(M_v(s)\). This doubles
 \(\log M_v(s)\) in \(b_v\), and the coefficient heights grow by a fixed
 factor.
\end{itemize}
Thus (iii) holds for \eqref{eq:comm}, with polynomial representatives of
degree at most four in \(\mathbf B\) and two in \(U\).
\end{remark}

It remains to verify (iv), non-triviality in the sense of
\cite[Definition~5.14]{QM}. We treat the single-chart systems in
\Cref{sec:single} and the paired systems in \Cref{sec:hecke,sec:frobenius}.

\section{Non-triviality of the single-chart systems}\label{sec:single}

\begin{proposition}\label{prop:single-nontrivial}
Each of the systems \eqref{eq:fixed}, \eqref{eq:zero} and
\eqref{eq:comm}, together with the chart equations, is non-trivial at
\(\eta(s)\) in the sense of \cite[Definition~5.14]{QM}.
\end{proposition}

\begin{proof}
Recall from \Cref{subsec:constant} that the coordinates vary, so the
families \(E_1,E_2,E_3\) are nonisotrivial and pairwise not isogenous at the
geometric generic point. Let \(\mathscr X\) be the single-chart variety,
on which \(t=x(Z)\), and let
\[
 p:\mathscr X\longrightarrow S_{\mathrm{aff}},\qquad
 (t,Z,P_e,P_f)\longmapsto Z,
\]
be the projection to the base. Recall that \cite[Definition~5.14]{QM}
quantifies over every irreducible closed subvariety
\(W\subseteq\mathscr X\) that contains \(\eta(s)\) and on which \(t\) is
nonconstant, and over every prime divisor \(\mathfrak p\) of its
normalization \(\widetilde W\) with
\(\operatorname{ord}_{\mathfrak p}(t)>0\).

Suppose that, for some \(W\) and \(\mathfrak p\), all series of the
system vanish in \(\widehat{\mathcal O}_{\widetilde W,\mathfrak p}\).
Then \(Z=\zeta(t)\) formally at \(\mathfrak p\). Hence, near \(\mathfrak p\), we have \(z(p(w))=\zeta(t(w))\) and
\(x(p(w))=t(w)\), and \(\mathcal T_i(t(w))\) is parallel
transport from \(E_{i,p(w)}\) to \(E_{i,s_0}\). By
\cite[Lemma~7.22]{QM}, these identities continue around loops on a
smooth dense open subset \(W^\circ\subseteq W\). The proof given there applies verbatim: part (i) uses the convergence of
the fixed series near \(t=0\) and standard facts on normal noetherian
domains, and part (ii) the analytic continuation of solutions of the
fixed differential systems; all of these are available here. The continuation of \(\mathcal T_i(t(w))\) along a path
\(\gamma\) in \(W^\circ\) is parallel transport along \(p\circ\gamma\). Since
\(x\circ p=t\) is nonconstant on \(W\), the morphism \(p|_W\) is
dominant. Each \(P_k\) used below is invertible at \(\eta(s)\), hence on a
dense open subset, and we choose the base point of the continuation
there. Its algebraic coordinates then return to their original values
around every loop.

Write \(V_i:=H_1(E_i,\mathbb Q)\) for the Betti local systems. The
identity component of the Zariski closure of the monodromy on the three factors is \(\prod_i\SL(V_i)\). This is shown in the proof of
\cite[Lemma~4.2]{DOP25}, and we recall the argument.
\begin{itemize}
\item For a single nonisotrivial family, the identity component is
 semisimple by Deligne's semisimplicity theorem. It is nontrivial, since
 finite monodromy would make the family isotrivial. Hence it is
 \(\SL(V_i)\).
\item For two families, consider a proper subgroup of
 \(\SL(V_i)\times\SL(V_j)\) projecting onto both factors. Up to the
 centers, it is the graph of an automorphism, which is inner. After a
 finite cover, \(V_i\cong V_j\) as local systems. By the theorem of the
 fixed part, \(\Hom(V_i,V_j)^{\pi_1}\) is then a one-dimensional rational
 Hodge structure of weight zero, hence of type \((0,0)\), and it contains
 an isogeny. This is a contradiction.
\item Goursat's lemma, for products of simple groups, then gives the full
 product.
\end{itemize}
The same argument applies over any smooth irreducible base on which the
families are nonisotrivial and pairwise not isogenous. By
\cite[Lemma~7.4]{QM}, applied to the dominant morphism
\(p:W^\circ\to S_{\mathrm{aff}}\), the identity component is unchanged by
pullback.

Consider first \eqref{eq:fixed} and \eqref{eq:zero}. For
\(k=(1,j)\in\{e,f\}\), the series \(F_k\) is a matrix coefficient of
\(V_j\otimes V_1^\vee\): it pairs the transported \(P_k\) with
\(B_k^\vee\) by the trace. Both arguments of the pairing are nonzero.
Since the translates of a nonzero vector span the irreducible module
\(V_j\otimes V_1^\vee\), the coefficient is not identically zero. Since
the central element \(-I\) of \(\SL(V_1)\) changes its sign, it is not
constant either. In particular it is identically equal neither to
\(T_k\) nor to zero.

Consider now \eqref{eq:comm}. We fix the first factor and let the
factors \(\SL(V_2)\) and \(\SL(V_3)\) act independently. Since the central
and moving maps are invertible, each \(M_k\) then ranges over all
matrices of its fixed nonzero determinant \(\delta_k\) \((k\in\{e,f\})\). In
particular \(M_e\) and \(M_f\) can take the values
\[
 \begin{pmatrix}\delta_e&1\\0&1\end{pmatrix}
 \qquad\text{and}\qquad
 \begin{pmatrix}\delta_f&0\\1&1\end{pmatrix},
\]
whose commutator has upper-left entry \(1\). This contradicts the
vanishing of all four entries after continuation.
\end{proof}

\section{Identities on the doubled chart}\label{sec:hecke}

Since the curve and the center are defined over a fixed number field,
only finitely many conjugates \(C^\tau\) occur. We
write \((j_1,j_2,j_3,j'_1,j'_2,j'_3)\) for the coordinates on
\(Y(1)^3\times Y(1)^3\).

\subsection{The two traces determine the matching of the three
coordinates}\label{subsec:matching}

Let \(\mathscr X_2\) be the doubled-chart variety, on which
\(t=x(Z)\) and \(t\mathcal R=x^\tau(Z^\tau)\), with the two projections
\[
 p:\mathscr X_2\to S_{\mathrm{aff}},\quad
 (t,U)\mapsto Z,\qquad
 p':\mathscr X_2\to S^\tau_{\mathrm{aff}},\quad
 (t,U)\mapsto Z^\tau .
\]

\begin{proposition}\label{prop:matching}
Let \(W\subseteq\mathscr X_2\) be an irreducible closed subvariety with
\(\eta(s)\in W\) and \(t|_W\) nonconstant, and let \(\mathfrak p\) be a
prime divisor of its normalization \(\widetilde W\) with
\(\operatorname{ord}_{\mathfrak p}(t)>0\). Suppose that both
differences \eqref{eq:paired} and all chart equations of the doubled
chart vanish in \(\widehat{\mathcal O}_{\widetilde W,\mathfrak p}\). Then
the closure of
the image of \(W\) in \(C\times C^\tau\) is a curve, dominant over
\(C\) and over \(C^\tau\), on which
\[
 \Phi_{n_1}(j_1,j'_1)=\Phi_{n_2}(j_2,j'_2)=\Phi_{n_3}(j_3,j'_3)=0
\]
for some \(n_1,n_2,n_3\ge1\). This curve contains the image of
\(\eta(s)\), that is, the pair of triples of \(j\)-invariants at \(s\) and
\(s^\tau\).
\end{proposition}

We call such a curve a \emph{dominant Hecke curve}.

\begin{proof}
The projection of \(W\) to the conjugate base is dominant as well.
Indeed, if \(t\mathcal R=x^\tau\circ p'\) were constant on \(W\), its
residue at \(\mathfrak p\) would make this constant zero, whereas its
value at \(\eta(s)\) is \(\tau(\xi)\ne0\). Denote by \(E_1,E_2,E_3\) and
\(E_1^\tau,E_2^\tau,E_3^\tau\) the two triples of elliptic families
pulled back to \(W\). By dominance, the three families in each triple
are pairwise not isogenous at the geometric generic point.

We group the six families by isogeny at the geometric generic point.
After a finite cover and restriction to a dense open subset, isogenies
within a class identify their rational local systems. For one
representative of each class, the connected monodromy is the product of
their \(\SL_2\)'s, by the argument recalled in the proof of
\Cref{prop:single-nontrivial}. When
\(\dim W>1\), we restrict to a smooth algebraic curve in this open subset
through a very general complex point, on which every coordinate
remains nonconstant. We choose the point outside the countable union of
modular loci for the pairs that are generically not isogenous on \(W\).
These pairs then remain generically not isogenous on the curve, and its
loops are loops in the open subset of \(W\), so the curve case applies.
Alternatively, the argument recalled in the proof of
\Cref{prop:single-nontrivial} applies directly over a smooth dense open
subset of \(W\).

Write \([E_i]\) and \([E_i^\tau]\) for the resulting isogeny classes.
For \(k=(1,j)\in\{e,f\}\), the central element \(-I\) in the factor of a
single class acts by \(-1\) on the trace coefficient \(F_k\) exactly
when the class is one of \([E_1],[E_j]\), and similarly for
\(F_k^\tau\). If a class occurred on just one side of the identity
\eqref{eq:paired}, this element would change the sign of just that side.
The identity on the connected monodromy group would then force that
coefficient to vanish identically. This contradicts the argument in the
proof of \Cref{prop:single-nontrivial}, which applies because \([E_1]\)
and \([E_j]\) are distinct factors. Applying this test to both pairs
\(e\) and \(f\), we obtain
\[
 \{[E_1],[E_2]\}=\{[E_1^\tau],[E_2^\tau]\}
 \qquad\text{and}\qquad
 \{[E_1],[E_3]\}=\{[E_1^\tau],[E_3^\tau]\}.
\]
Intersecting gives \([E_1]=[E_1^\tau]\), and hence \([E_2]=[E_2^\tau]\)
and \([E_3]=[E_3^\tau]\). Thus each \(E_i\) is geometrically generically
isogenous to \(E_i^\tau\), with the same label \(i\).

In each coordinate we choose a generic cyclic isogeny, of degree
\(n_i\). The three modular equations \(\Phi_{n_i}(j_i,j'_i)=0\) extend
over the closure of the image by closedness. Their common locus in
\(Y(1)^3\times Y(1)^3\) is finite over either factor, so the closure of
the image is a curve, dominant over \(C\) and \(C^\tau\). It contains
the image of \(\eta(s)\), since \(\eta(s)\in W\).
\end{proof}

\subsection{Uniform degrees for the coordinatewise isogenies}\label{subsec:degrees}

\begin{proposition}\label{prop:uniform-degrees}
There is a finite set \(\mathcal M\subset\mathbb Z_{>0}^{\,3}\),
depending only on \(C\) and its conjugates,
with the following property. At every point of every dominant Hecke
curve in \(C\times C^\tau\), for every \(\tau\), there are isogenies
in the three coordinates of degrees \((m_1,m_2,m_3)\in\mathcal M\).
\end{proposition}

The proof adapts \cite[Lemma~7.13 and Proposition~7.10]{QM}.

\begin{proof}
Pass to a fixed torsion-free level, and take one pair of irreducible
lifted algebraic curves. Choose irreducible analytic components
\(\mathcal C,\mathcal C'\subset\mathfrak H^3\) of their preimages. Let
\(\Gamma\) be the level group, and \(H_{\mathcal C}\) the stabilizer of
\(\mathcal C\) in \(\PSL_2(\mathbb R)^3\).

We first show that \(H_{\mathcal C}\) is discrete. The closedness
argument of \cite[Lemma~7.13]{QM} shows that \(H_{\mathcal C}\) is a
closed Lie subgroup. The projective monodromy lies in
\(H_{\mathcal C}\cap\Gamma\) and is Zariski dense in the product, so
\(\Lie H_{\mathcal C}\) is an ideal of \(\mathfrak{sl}_2(\mathbb R)^3\).
A nonzero ideal contains a whole simple factor. Its orbit through a
smooth point would give a one-dimensional coordinate disc in
\(\mathcal C\) on which the other two coordinates are constant. By
irreducibility they would then be constant on \(\mathcal C\), contrary to the assumption that the three coordinates vary.

Next, \([H_{\mathcal C}:H_{\mathcal C}\cap\Gamma]<\infty\). As in the
finite-index part of \cite[Lemma~7.13]{QM}, \(H_{\mathcal C}\) acts
properly on the normalization of \(\mathcal C\), with finite point
stabilizers. The quotient by \(H_{\mathcal C}\cap\Gamma\) is the
normalization of the algebraic curve at level, which has finite intrinsic
Poincar\'e area. A small disc and its translates under
\(H_{\mathcal C}\) give disjoint discs of equal positive area in this
quotient, indexed by the double cosets of \(H_{\mathcal C}\cap\Gamma\)
and a finite point stabilizer in \(H_{\mathcal C}\). Hence the index is
finite.

Finally, consider a dominant Hecke curve. On a simply connected open
subset of an isogeny cover, the two lifted triples of period maps
\(\varpi,\varpi'\) satisfy \(\varpi'=g\varpi\) for a fixed rational tuple
\(g\in\PGL_2^+(\mathbb Q)^3\). Continuation gives
\(g\mathcal C_1=\mathcal C'_1\) for irreducible components
\(\mathcal C_1,\mathcal C'_1\) of the preimages of the two curves at
level. This is the graph argument of \cite[Proposition~7.10]{QM}, applied
in each coordinate.

Replacing \(g\) by \(\gamma'g\gamma^{-1}\), with \(\gamma,\gamma'\in\Gamma\),
does not change the image of the graph in the product of modular curves.
Moreover, \(\Gamma\) permutes transitively the components above each
irreducible curve at level. We may therefore take
\((\mathcal C_1,\mathcal C'_1)=(\mathcal C,\mathcal C')\) from a fixed
finite list of pairs.

If a rational tuple \(g_0\) takes \(\mathcal C\) to \(\mathcal C'\), then
every tuple doing so lies in \(g_0H_{\mathcal C}\), and right
multiplication by \(H_{\mathcal C}\cap\Gamma\) preserves the image of the
graph. Thus only the finitely many cosets in
\(g_0H_{\mathcal C}/(H_{\mathcal C}\cap\Gamma)\) are relevant, for each
of the finitely many pairs of lifts. We choose a rational representative
of each and clear denominators in its three homology matrices; the
resulting positive determinants form \(\mathcal M\). The corresponding
isogenies exist throughout each graph, and the Hecke relations of fixed
degree are closed, so they exist at every point of its algebraic
closure.
\end{proof}

We note that \Cref{prop:uniform-degrees} asserts the existence of one
triple of maps at each point, not a bound for all of them.

\section{Frobenius on the two reduced endomorphisms}\label{sec:frobenius}

We retain the notation of \Cref{sec:hecke}, and suppose that the pair of
triples at \(s\) and \(s^\tau\) lies on a dominant Hecke curve. By
\Cref{prop:uniform-degrees}, we may choose isogenies
\[
 \phi_{i,s}:E_{i,s}\longrightarrow E_{i,s}^{\tau},\qquad
 \deg\phi_{i,s}=m_i\quad(i=1,2,3),\qquad (m_1,m_2,m_3)\in\mathcal M,
\]
where \(E_{i,s}^\tau:=\tau(E_{i,s})\) is the fiber of \(E_i^\tau\) at
\(s^\tau\). Fix a good embedding \(\rho\) in the class of the paired
system labelled by \(\tau\), and take geometric reductions at one
extension of \(\rho\). We extend the maps over the traits of good
reduction.

Since \(\rho\tau=\sigma_\rho\rho\) on \(L_s\), we have
\(\rho(s^\tau)=\sigma_\rho(\rho(s))\). The automorphism \(\sigma_\rho\) of
\(\Qbar_p\) is continuous and induces \(u\mapsto u^p\) on residue fields.
Hence it maps the residue disc of \(\rho(s_0)\) onto that of
\(\rho(s_0^\tau)\), and the reduction of \(\rho(\tau(X))\) is the Frobenius
twist \((\widetilde X)^{(p)}\) of the reduction of \(\rho(X)\), for every
\(X\) defined over \(L_s\). Consequently, the selected discs identify the
reductions at \(s\) and \(s^\tau\) with those at \(s_0\) and \(s_0^\tau\),
respectively, and the reduction of \(E^\tau_{i,s_0^\tau}\) is
\(\widetilde E_i^{(p)}\). Thus the reduction \(\bar\phi_i\) of
\(\phi_{i,s}\) is a specified isogeny
\(\widetilde E_i\to\widetilde E_i^{(p)}\) between the reductions of two
fixed central curves.

\subsection{Lifting the specified maps}\label{subsec:lifting}
The following lemma is the elliptic analogue of \cite[Lemma~8.13]{QM}.
Its excluded primes are those placed in \(\Sigma\) in
\Cref{subsec:sigma}.

\begin{lemma}\label{lem:lifting}
There exist a finite extension of the field of definition of the central
data and a finite set of rational primes, depending only on the fixed
central curves, their finitely many conjugates and \(\mathcal M\), with
the following property. The set contains \(2\), \(3\), the primes
dividing an entry of an element of \(\mathcal M\), and the primes of bad
reduction of the central curves and their conjugates. At every embedding above a prime outside this
set, each \(\bar\phi_i\) is the reduction of an isogeny
\(\phi_{i,0}:E_{i,s_0}\to E_{i,s_0}^{\tau}\) of degree \(m_i\).
\end{lemma}

\begin{proof}
Fix an ordered pair \((E,E')\) of central elliptic curves from the fixed
list, and an integer \(\mu\) occurring as an entry of an element of
\(\mathcal M\). Away from the primes dividing \(\mu\), every possible
kernel of an isogeny of degree \(\mu\) in the geometric special fiber
lifts, after a fixed extension, to a subgroup \(J\subset E[\mu]\) of order
\(\mu\). For each of the finitely many quotients \(E/J\), we exclude the
primes at which \(j(E/J)\ne j(E')\) become equal. For the quotients
isomorphic to \(E'\), we choose an isomorphism and all geometric
automorphisms over a fixed extension. Away from \(2,3\) and the finitely
many primes at which \(j(E')\) specializes to \(0\) or \(1728\), the
geometric automorphisms of \(E'\) specialize onto those of its
reduction. Hence the prescribed isomorphism of the reduced quotient with
\(\overline{E'}\) lifts, and gives the prescribed map. We take the union
of these exclusions, the bad primes and the primes dividing any \(m_i\).
All lifting choices may be made over algebraic closures at the embedding
in question.
\end{proof}

\subsection{The two signs}\label{subsec:signs}
Write \(E_{i,0}:=E_{i,s_0}\), and for \(k\in\{e,f\}\) set
\(a_{k,s}:=a_k\) and \(a_{k,0}:=b_k\). Let \(k=(1,j)\in\{e,f\}\) and
\(r\in\{s,0\}\). By the non-CM hypothesis,
\(\Hom^0(E_{1,r},E_{j,r})=\mathbb Q\,a_{k,r}\). Since conjugation
preserves degrees, there is \(c_{k,r}\in\mathbb Q^\times\) such that
\begin{equation}\label{eq:rank-one}
 \phi_{j,r}^{-1}\,a_{k,r}^\tau\,\phi_{1,r}=c_{k,r}\,a_{k,r}
 \qquad\text{and}\qquad c_{k,r}^2=m_1/m_j,
\end{equation}
where the second identity follows on taking degrees.

Let \(\Fr_i:\widetilde E_i\to\widetilde E_i^{(p)}\) be the relative
\(p\)-Frobenius. As shown at the beginning of this section, the
reduction of the conjugate fiber is \(\widetilde E_i^{(p)}\). We set
\[
 \theta_i:=\bar\phi_i^{-1}\Fr_i\in\End^0(\widetilde E_i)^\times .
\]

\begin{lemma}\label{lem:signs}
For \(k\in\{e,f\}\) we have
\begin{equation}\label{eq:signs}
 \theta_1D_k\theta_1^{-1}=\epsilon_kD_k,\qquad\text{where}\quad
 \epsilon_k:=c_{k,s}/c_{k,0}\in\{1,-1\}.
\end{equation}
\end{lemma}

\begin{proof}
Write \(k=(1,j)\). Reduction of homomorphisms is injective and
compatible with composition and degrees. After tensoring with
\(\mathbb Q\), it is therefore compatible with inverses in \(\Hom^0\).
Naturality of Frobenius gives
\(\Fr_j\,\bar a_{k,r}\,\Fr_1^{-1}=\overline{a_{k,r}^\tau}\). Reducing
\eqref{eq:rank-one}, with the same \(\bar\phi_i\) for \(r=s\) and \(r=0\)
by \Cref{lem:lifting}, we obtain
\[
 \theta_j\,\bar a_{k,r}\,\theta_1^{-1}=c_{k,r}\,\bar a_{k,r}
 \qquad(r=s,0).
\]
Inverting the identity for \(r=0\), whose isogeny is \(b_k\), and
multiplying it by the identity for \(r=s\) gives \eqref{eq:signs}. Since
\(c_{k,s}^2=c_{k,0}^2\), the ratio \(\epsilon_k\) is a sign.
\end{proof}

\begin{corollary}\label{cor:reassigned}
The embedding \(\rho\) was assigned to the commutator system
\eqref{eq:comm} or to a system \eqref{eq:zero}.
\end{corollary}

\begin{proof}
If \(\epsilon_k=-1\) for some \(k\in\{e,f\}\), invariance of the reduced
trace under conjugation gives \(\tr D_k=0\), and hence
\(F_k(\xi,A_k)=N_{0,k}\tr D_k=0\) at this embedding.

If both signs are positive, both \(D_e\) and \(D_f\) centralize
\(\theta_1\). The
degree \(p/m_1\) of \(\theta_1\) has odd \(p\)-adic valuation, since
\(p\nmid m_1\), so \(\theta_1\) is not a rational scalar. If the reduction
is ordinary, the endomorphism algebra is an imaginary quadratic field
\cite[\S\S III.9 and~V.3]{Sil09}.
If it is supersingular, the centralizer of the nonscalar \(\theta_1\) in
the quaternion algebra is the quadratic field \(\mathbb Q(\theta_1)\). In
either case the centralizer is commutative, so that
\(D_eD_f=D_fD_e\).

In both cases the partition of \Cref{subsec:partition} assigned \(\rho\)
before reaching the paired system.
\end{proof}

\subsection{Proofs of the main results}\label{subsec:proofs}

\begin{proof}[Proof of \Cref{thm:height}]
If a coordinate is constant, the bound follows from
\Cref{rem:constant-coordinates}. Otherwise, by \Cref{subsec:variables} we may
assume \(s\in\mathcal S^\circ\). If
\(h<Q_{\mathrm{thr}}(d)\) there is nothing to prove, so we may assume
the contrary, and let \(\mathcal V\) be the class of
\Cref{prop:proximity}. By \Cref{subsec:hypotheses}, it satisfies
hypotheses (i)--(iii) of \cite[Theorem~5.17]{QM}. If it is a single-chart
class, (iv) holds by \Cref{prop:single-nontrivial}. Suppose that it is a
nonempty class of the paired system that fails
\cite[Definition~5.14]{QM}. Then \Cref{prop:matching} puts the pair of
triples at \(s\) and \(s^\tau\) on a dominant Hecke curve, and by
\Cref{cor:reassigned} every embedding of the class was assigned to an
earlier class, which is a contradiction. Hence (iv) holds in all cases,
and \cite[Theorem~5.17]{QM} gives \(h_S(s)\ll d^c\).
\end{proof}

\begin{proof}[Proof of \Cref{thm:two-relations}]
By \Cref{thm:height} and the proof of
\cite[Proposition~5.15]{DO25}, the larger of the two minimal cyclic
isogeny degrees at \(s\in\mathcal S\) is bounded by a power of \(d\).
This lower bound for Galois orbits replaces \cite[Lemma~4.2]{HP12} in the
counting argument of \cite[Proposition~5.1]{HP12}, as in \cite[\S1.A]{DO25} and \cite[\S5]{DOP25}. That proposition
treats one pair of index pairs at a time. Its proof uses the asymmetry
hypothesis only through \cite[Lemma~4.2]{HP12}, and it works with any
positive exponent in place of \(1/6\); compare the remark preceding
\cite[Lemma~4.2]{HP12}. The fixed finite cover \(S\to C\) changes the
degrees of points by a bounded factor, so the counting argument gives finiteness of the points in question
whose coordinates are not CM. Together with the reduction to non-CM
points in \Cref{subsec:setup}, this proves the theorem.
\end{proof}

\begin{proof}[Proof of \Cref{thm:main}]
If \(C\) is not defined over \(\Qbar\), this is \cite[Theorem~1.4]{Pil17}.
Suppose that \(C\) is defined over \(\Qbar\). A special subvariety of
codimension at least two is contained in one of codimension two, defined
by two conditions. Each condition asserts either that a coordinate is a
singular modulus or that two coordinates are related by a modular
polynomial; see \cite[\S5.4]{HP12}.

If one of the two conditions is of the first kind, the subvariety is not
strongly special. By \cite[Theorem~2]{HP12}, only finitely many points
of \(C\) lie on such subvarieties.

Suppose that both conditions are modular relations
\(\Phi_M(x_i,x_j)=0\) and \(\Phi_N(x_{i'},x_{j'})=0\). If
\(\{i,j\}=\{i',j'\}\), then \(M\ne N\), since the subvariety has
codimension two, and both coordinates \(x_i,x_j\) are special on each of
its components. Indeed, if the elliptic curves with \(j\)-invariants
\(x_i\) and \(x_j\) were isogenous and not CM, their group of
homomorphisms would be generated by one isogeny \(\psi\). A cyclic
isogeny \(m\psi\) has \(|m|=1\), since its kernel contains the
non-cyclic group \(E[m]\) otherwise, so all cyclic isogenies between
them would have degree \(\deg\psi\). These components are therefore
covered by the previous case.

Otherwise \(\{i,j\}\ne\{i',j'\}\), and the two pairs share an index,
since there are only three coordinates. After permuting the coordinates, which preserves the hypothesis
on \(C\), the conditions become \(\Phi_M(x_1,x_2)=\Phi_N(x_1,x_3)=0\),
and \Cref{thm:two-relations} applies. There are finitely many choices
of index pairs.
\end{proof}

\begin{remark}\label{rem:three-coordinates}
\Cref{thm:two-relations} extends to curves in \(Y(1)^n\). Let
\(C\subset Y(1)^n\), with \(n\ge3\), be an irreducible curve over
\(\mathbb C\) that is not contained in any proper special subvariety, and
let \(i_1,i_2,i_3\) be distinct indices. Then only finitely many points of
\(C\) satisfy
\[
 \Phi_M(x_{i_1},x_{i_2})=\Phi_N(x_{i_1},x_{i_3})=0
 \qquad\text{for some }M,N\ge1.
\]
Indeed, let \(C_I\subset Y(1)^3\) be the closure of the projection of
\(C\) to the coordinates \(I=\{i_1,i_2,i_3\}\).
\begin{itemize}
\item If \(C_I\) is a point, no point of \(C\) satisfies the two
 relations: otherwise they would hold on all of \(C\), and \(C\) would lie
 in a proper special subvariety.
\item We may therefore assume that \(C_I\) is a curve. It is not
 contained in a proper special
 subvariety of \(Y(1)^3\), since the preimage of such a subvariety is a
 proper special subvariety of \(Y(1)^n\) containing \(C\).
\item \Cref{thm:main} applies to \(C_I\), and the projection
 \(C\to C_I\) has finite fibers.
\end{itemize}
\end{remark}

\subsection*{Acknowledgments}
The author heartily thanks Gal Binyamini for interesting discussions and
his encouragement. The author received funding from the European Union
(ERC, SharpOS, 101087910) and from the Israel Science Foundation (grant
No.~2067/23). Several of the ideas here originated while the author was
at the Institute for Advanced Study, supported by the Minerva Research
Foundation Member Fund for the academic year 2025--26.

\subsection*{Declaration of AI usage}
While preparing the notes \cite{QM} on the papers released by OpenAI, I
realized that the techniques developed therein, combined with the ideas
of \cite{DOP25} and two previously unpublished observations, were
sufficient to establish \Cref{thm:main}.

On November 18, 2026, I described the proof strategy to ChatGPT,
including the two crucial observations
(\Cref{prop:matching,cor:reassigned}), which date back to May 2026 and
build on the arguments in \cite{DOP25}. I asked ChatGPT to verify that
the lemmas established in the notes were sufficient to carry out this
strategy. I subsequently used Claude to perform an initial clean-up of
the manuscript.

Both the present manuscript and the notes \cite{QM} remain in
preliminary form, and I intend to undertake substantial revisions of
both in the coming weeks.

Following the advice of senior mathematicians in my field, I prepared a
preliminary manuscript with the intention of making it publicly
available should a claimed proof of the same result appear online.

\end{document}